\documentclass[titlepage,twoside,12pt]{article}
\usepackage{amssymb}
\usepackage{amsfonts}
\usepackage{amsmath}
\begin{document}

{\bf \Large The Iterative Simplicity of Basic \\ \\ Topological Operations} \\ \\

{\bf Elemer E Rosinger} \\
Department of Mathematics \\
and Applied Mathematics \\
University of Pretoria \\
Pretoria \\
0002 South Africa \\
eerosinger@hotmail.com \\ \\

{\bf Abstract} \\

Semigroups generated by topological operations such as closure, interior or boundary are
considered. It is noted that some of these semigroups are in general finite and
noncommutative. The problem is formulated whether they are always finite. \\ \\

{\bf 1. Iterating Closure, Interior and other Topological \\
        \hspace*{0.4cm} Operations} \\

Let $X$ be a topological space with the set ${\cal T} \subseteq {\cal P} ( X )$ of open
subsets. For convenience, we shall denote by $c ( A )$ and $i ( A )$ the closure,
respectively, interior of a subset $A \subseteq X$. \\

As is well known, some of the iterates of $c$ and $i$ play an important role in topology. For
instance, a subset $A \subseteq X$ is called {\it nowhere dense}, iff $i ( c ( A ) ) = \phi$.
Further, a countable union of nowhere dense subsets is called of {\it first Baire category},
the essential fact in this regard being that no complete metric space is of first Baire
category. \\

Motivated by the above, here various iterates of $c$ and $i$, as well as of other related
basic operations will be considered. In this regard we can note that all such operations are
mappings of ${\cal P} ( X )$ into itself, thus their compositions are associative.
Consequently, we can consider the {\it free semigroup} {\bf \large S} $( {\cal T} )$ of
mappings of ${\cal P} ( X )$ into itself generated by these topological operations which will
be listed below. This semigroup is in general obviously {\it noncommutative}. \\

The nontrivial aspect involved is that, in view of well known relations, such as \\

(1.1)~~~ $ c^2 = c,~~~ i^2 = i $ \\

a number of elements in this noncommutative semigroup {\bf \large S} $( {\cal T} )$ correspond
to the same mappings of ${\cal P} ( X )$ into itself. Therefore, their identification, that is,
the identification of the different elements in {\bf \large S} $( {\cal T} )$ is of
interest. \\
Here we shall be interested in such identification which hold for {\it all} topological spaces
$( X, {\cal T})$. \\
Obviously, in the case of particular topological spaces, one can find more such
identifications. \\

Let us list now the other mappings of ${\cal P} ( X )$ which we shall consider. \\

One of them is the operation of taking the complementary, namely, $\complement ( A ) = X
\setminus A$, with $A \subseteq X$. Here we have further reductions in the different elements
of {\bf \large S} $( {\cal T} )$, since \\

(1.2)~~~ $ c ( \complement ( A ) ) = \complement ( i ( A ) ),~~~
                      i ( \complement ( A ) ) = \complement ( c ( A ) ),~~~ A \subseteq X $ \\

hence \\

(1.3)~~~ $ c \, \complement = \complement \, i,~~~ i \, \complement = \complement \, c $ \\

One also defines the exterior $e ( A )$ of a subset $A \subseteq X$, given by \\

(1.4)~~~ $ e = i \, \complement $ \\

Another important topological operation is that of boundary of a subset $A \subseteq X$, which
we shall denote by $b ( A )$. Here we have the obvious relation \\

(1.5)~~~ $ b ( A ) = b ( \complement ( A ) ),~~~ A \subseteq X $ \\

or simply \\

(1.6)~~~ $ b \, \complement = b $ \\

Also \\

(1.7)~~~ $ b ( A ) =
            \complement ( i ( A ) \cup i ( \complement  ( A ) ) ),~~~ A \subseteq X $ \\

thus \\

(1.8)~~~ $ b = \complement ( i \cup i \, \complement ) $ \\

or in view of (1.4) \\

(1.7)~~~ $ b = \complement ( i \cup e ) $ \\

Let us recall that a subset $A \subseteq X$ is called a {\it boundary set}, iff \\

(1.8)~~~ $ i ( A ) = \phi $ \\

thus $A$ is nowhere dense, iff \\

(1.9)~~~ $ i ( c ( A ) ) = \phi $ \\

Let us now consider the related operations \\

(1.10)~~~ $ b_i ( A ) =  b ( A ) \cap A,~~
b_e ( A ) = b ( A ) \cap \complement ( A ),~~~ A \subseteq X $ \\

called respectively the {\it internal} and {\it external} boundary of $A$. Clearly, we have \\

(1.11)~~~ $ b = b_i \cup b_e $ \\

There are several further frequently used topological operations. One of them is the {\it
derived} set $d ( A )$ of a subset $A \subseteq X$, defined by \\

(1.12)~~~ $ d ( A ) = \{ x \in X ~|~
            \forall~~ G \in {\cal T} ~:~ ( A \cap G ) \setminus \{ x \} \neq \phi ~\} $ \\

And now we define \\

(1.13)~~~ {\bf \large S} $( {\cal T} ) $ \\

as the free semigroup generated by the set of mappings $\{ c, i, \complement, e, b, b_i, b_e,
d \}$, and as customary with semigroups, we assume that it contains the neutral element $id_X$
which maps each $A \subseteq X$ into itself. \\ \\

{\bf 2. A Simpler Problem} \\

Let us start with the simpler problem of studying the subsemigroup of {\bf \large S}
$( {\cal T} )$ which is generated by the two operations $c$ and $i$ alone. This subsemigroup,
in view of (1.1), is obviously given as follows \\

(2.1)~~~ {\bf \large S}$_{c,\, i} ( {\cal T} ) =
                            \{~ id_X, c, i, ci, ic, ici, cic, cici, icic, \ldots ~\} $ \\

Let us consider the partial order relation $\alpha \rightarrow \beta$ between mappings of
${\cal P} ( X )$ into itself, defined by \\

(2.2)~~~ $ \alpha ( A ) \subseteq \beta ( A ),~~~ A \subseteq X $ \\

Then obviously \\

(2.3)~~~ $ i \rightarrow id_X \rightarrow c $ \\

On the other hand, we have \\

(2.4)~~~ $ i ( A ) \subseteq i ( B ),~~~
              c ( A ) \subseteq c ( B ),~~~ A \subseteq B \subseteq X $ \\

{\bf Lemma 2.1.} \\

{\bf \large S}$_{c,\,i} ( {\cal T} )$ is a partially ordered semigroup with respect to
$\rightarrow$, and in general it is noncommutative. \\

{\bf Proof.} \\

Let $\alpha, \beta, \gamma \in$ {\bf \large S}$_{c,\,i} ( {\cal T} )$, with $\alpha
\rightarrow \beta$. \\
We show that $\gamma \, \alpha \rightarrow \gamma \, \beta$. Let $A \subseteq X$. Then (2.2)
gives $\alpha ( A ) \subseteq \beta ( A )$, hence (2.1), (2.4) result in $\gamma ( \alpha
( A ) ) \subseteq \gamma ( \beta ( A ) )$. \\
Similarly,  $\alpha \, \gamma \rightarrow \beta \, \gamma$. Indeed, for  $A \subseteq X$, we
have $\gamma ( A ) \subseteq X$, hence (2.2) gives $\alpha ( \gamma ( A ) ) \subseteq \beta
( \gamma ( A ) )$. \\

We show that, in general, none of the relations holds \\

(2.5)~~~ $ i\, c = c\, i,~~~ i\, c \rightarrow c\, i,~~~ c\, i \rightarrow i\, c $ \\

Let $( X, {\cal T}) = \mathbb{R}$. If $A = [ 0, 1 ] \cap \mathbb{Q}$, then $c ( i ( A ) ) =
\phi \subsetneqq ( 0, 1 ) = i ( c ( A ) )$. Thus the first and the second of the above
relations do not hold. Let now $A = [ 0, 1 ]$, then $i ( c ( A ) ) = ( 0, 1 ) \subsetneqq
[ 0, 1 ] = c ( i ( A ) )$, hence the third relation above cannot hold. \\

{\bf Lemma 2.2} \\

The relations hold \\

(2.6)~~~ $ \begin{array}{l}
                i \rightarrow c i \rightarrow c \\
                i \rightarrow i c \rightarrow c \\
                i \rightarrow i c i \rightarrow i c \rightarrow c \\
                i \rightarrow c i \rightarrow c i c \rightarrow c \\
                i \rightarrow c i \rightarrow c i c i \rightarrow c i c \rightarrow c \\
                i \rightarrow i c i \rightarrow i c i c \rightarrow i c \rightarrow c \\
             \end{array} $ \\

Consequently we have \\

(2.7)~~~ $ c i c i = c i,~~~ i c i c = i c $ \\

{\bf Proof.} \\

For the first relation in (2.6) we compose (2.3) on the left with $c$ and obtain $c i
\rightarrow c \rightarrow c^2$, thus $c i \rightarrow c $. Composing now (2.3) on the right
with $i$, the result is $i^2 \rightarrow i \rightarrow c i$, or $i \rightarrow c i$. \\

The second relation in (2.6) follows by composing (2.3) on the right with $c$, and thus
obtaining $i c \rightarrow c \rightarrow c^2$ or $i c \rightarrow c$. While composing
(2.3) on the left with $i$, it follows that $i^2 \rightarrow i \rightarrow i c$, or $i
\rightarrow i c$. \\

For the third relation in (2.6) we compose the first relation in it with $i$ on the left and
obtain $i^2 \rightarrow i c i \rightarrow i c$, and then recall the second relation in
(2.6). \\

The fourth relation in (2.6) is obtained by composing the second relation in it with $c$ on
the left, with the result $c i \rightarrow c i c \rightarrow c^2$, and then use the first
relation in (2.6).\\

The fifth relation in (2.6) results from the composition on the left with $c$ of the third
relation in (2.6), which gives $c i \rightarrow c i c i \rightarrow c i c \rightarrow c^2$,
after which we recall the first relation in (2.6). \\

The sixth relation in (2.6) comes from composing the fourth relation in (2.6) with $i$ on
the left, thus having $i^2 \rightarrow i c i \rightarrow i c i c \rightarrow i c$, and then
recalling the second relation in (2.6). \\

For (2.7) we proceed as follows. \\

We have $c i \rightarrow c i c i$ from the fifth relation in (2.6). On the other hand, the
fourth relation in (2.6) gives $c i c \rightarrow c$, which composed on the right with $i$,
yields $c i c i \rightarrow c i$. \\

As for $i c i c \rightarrow i c$, it follows from the sixth relation in (2.6). Now the third
relation in (2.6) gives $i \rightarrow i c i$, which composed with $c$ on the right results
in $i c \rightarrow i c i c$.

\hfill $\Box$ \\

The above lead now to \\

{\bf Theorem 2.1.} \\

The typically noncommutative semigroup generated by the closure and interior operations $c$,
respectively $i$, is given by \\

(2.8)~~~ {\bf \large S}$_{c,\, i} ( {\cal T} ) =
                            \{~ id_X, c, i, ci, ic, ici, cic ~\} $ \\

thus it has at most seven elements. With respect to the partial order relation $\rightarrow$,
these elements are in general arranged as follows \\

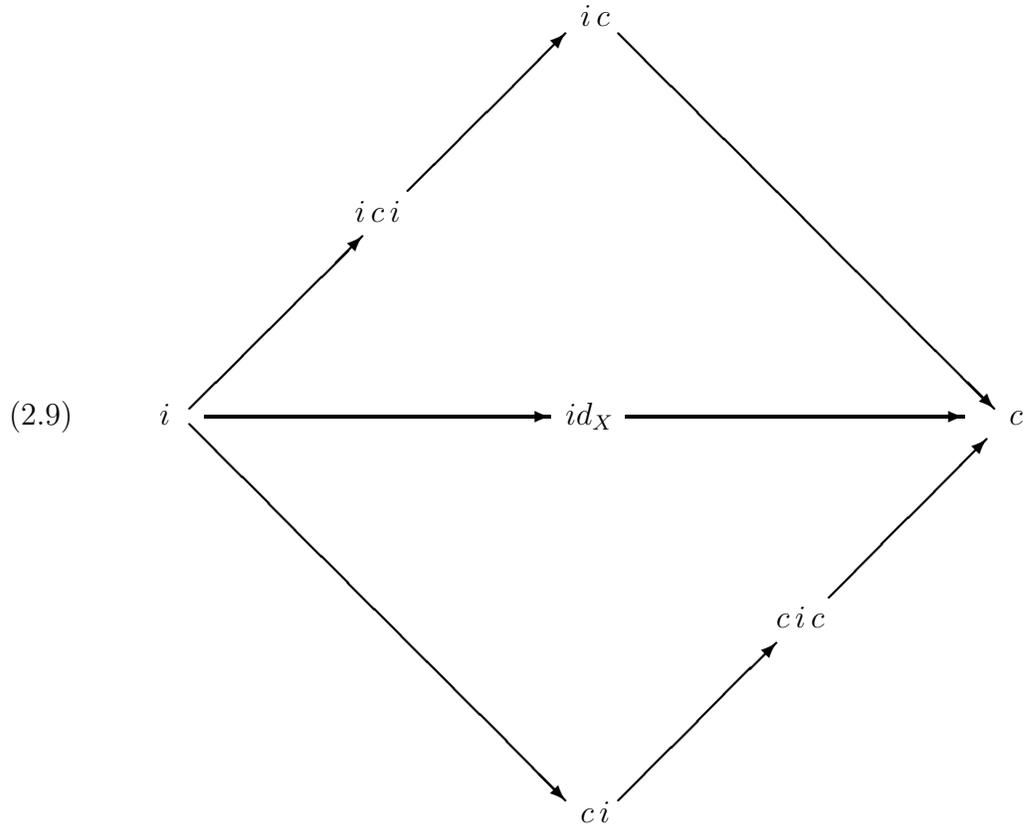
\begin{math}
\setlength{\unitlength}{1cm}
\thicklines
\begin{picture}(30,12)

\put(0,6){$(2.9)$}
\put(2,6){$i$}
\put(2.4,6.2){\vector(1,1){2.3}}
\put(4.6,8.7){$i\,c\,i$}
\put(5.3,9.1){\vector(1,1){2.1}}
\put(2.4,6){\vector(1,-1){5}}
\put(7.6,11.3){$i\,c$}
\put(7.6,0.7){$c\,i$}
\put(8.1,11.2){\vector(1,-1){5}}
\put(8.1,1){\vector(1,1){2.1}}
\put(10.2,3.3){$c\,i\,c$}
\put(10.9,3.7){\vector(1,1){2.1}}
\put(13.3,6){$c$}
\put(2.6,6.1){\vector(1,0){4.6}}
\put(7.4,6){$id_X$}
\put(8.2,6.1){\vector(1,0){4.5}}

\end{picture}
\end{math}

{\bf Remark 2.1.} \\

We typically have the relations \\

(2.10)~~~ $ c i c \neq c i,~~~ i c i \neq i c $ \\

as the following simple example shows it. Let $( X, {\cal T}) = \mathbb{R}$. If $A = [ 0, 1 ]
\cap \mathbb{Q}$, then $c ( i ( A ) ) = \phi$, while $c ( i ( c ( A ) ) ) = [ 0, 1 ]$. Also
$i ( c ( A ) ) = ( 0, 1 )$, while $i ( c ( i ( A ) ) ) = \phi$. \\ \\

{\bf 3. A Larger Semigroup} \\

Let us consider now the subsemigroup {\bf \large S}$_{c,\, i,\, \complement} ( {\cal T} )$
of {\bf \large S} $( {\cal T} )$ which is generated by the three operations $c$, $i$ and
$\complement$. This subsemigroup, in view of (1.3), (2.8) and of the relation $\complement \,
\complement = id_X$, is obviously given by \\

(3.1)~~~ {\bf \large S}$_{c,\, i,\, \complement} ( {\cal T} ) =
                \{~ id_X,  c, i, ci, ic, ici, cic, \complement, \complement c, \complement i,
                  \complement ci, \complement ic, \complement ici, \complement cic ~\} $ \\

thus in general, it has at most fourteen elements. \\ \\

{\bf 4. An Open Problem} \\

A next issue is the structure of the generally noncommutative semigroup \\

(4.1)~~~ {\bf \large S}$_{c,\, i,\, \complement, b} ( {\cal T} ) $ \\

generated by the five operations $c,\, i,\, \complement, b$. In view of relations such as
(1.6) and (1.8), one can expect a certain simplification in the structure of that semigroup.
This leads, among others, to the following \\

{\bf Open Problem} \\

1) Is the noncommutative semigroup {\bf \large S}$_{c,\, i,\, \complement, b} ( {\cal T} )$
finite, in the case of general topological spaces $( X, {\cal T} )$ ? \\

2) The same question for the noncommutative semigroup {\bf \large S} $( {\cal T} )$ in
(1.13). \\ \\

\end{document}